\documentclass[onecolumn]{autart}
\usepackage{graphicx}
\input{amssym}

\newcommand{\p}{\rm \partial}

\renewcommand{\theequation}{\arabic{section}.\arabic{equation}}
\begin{document}
\begin{frontmatter}
\title{Classical and Nonclassical symmetries of the\\ (2+1)- dimensional Kuramoto-Sivashinsky equation} 
\thanks[footnoteinfo]{ Corresponding author: Tel. +9821-73913426.
Fax +9821-77240472.}
\author[]{Mehdi Nadjafikhah\thanksref{footnoteinfo}}\ead{m\_nadjafikhah@iust.ac.ir},
\author[]{Fatemeh Ahangari}\ead{fa\_ahanagari@iust.ac.ir}

\address{School of Mathematics, Iran University of Science and Technology, Narmak, Tehran 1684613114, Iran.}
\renewcommand{\sectionmark}[1]{}
\begin{abstract}
%
%
In this paper,  we have studied the problem of determining the
largest possible set of symmetries for an important example of
nonlinear dynamical system: the Kuramoto-Sivashinsky (K-S) model
in two spatial and one temporal dimensions. By applying the
classical symmetry method for the K-S model, we have found the
classical symmetry operators. Also, the structure of the Lie
algebra of symmetries is discussed and the optimal system of
subalgebras of the equation is constructed. The Lie invariants
associated to the symmetry generators  as well as the
corresponding similarity reduced equations are also pointed out.
By applying the nonclassical symmetry method for the  K-S model
we concluded that the analyzed model do not admit supplementary,
nonclassical type, symmetries. Using this procedure, the
classical Lie operators only were generated.
\end{abstract}
\end{frontmatter}
\section{Introduction}
As it is well known, commonly the growth process of thin films is
from equilibrium and involves the interactions of a large number
of particles. Notable attentions  both in theory and experiment
have been devoted to the  study of the mechanism of thin film
growth (see \cite{[1]} and references therein). In order to
describe the growth process mathematically with the partial
differential equations of the surface height with respect to time
and spatial cooordinates, some models have been proposed by
researches. These models are based on the collective and
statistic behaviour of deposited particles. One of them include
the Kuramoto-Sivashinsky (K-S) models in which the surface
tension, diffusion, nonlinear effect and random deposition are
considered
and so on.\\
The K-S equation has been derived e.g. in the context of chemical
turbulence \cite{[2]}. As we declarde above, it is also so
important because it can describe the flow of a falling fluid
film \cite{[3]}.  The solutions of  the K-S equation can clarify
the features  about the evolution of surface morphology  which is
instructive to the  investigations of the mechanism of thin film
growth.\\
The two dimensional  Kuramoto-Sivashinsky (K-S)
 model represents a nonlinear dynamical system  which is defined in a two
dimensional space $\{x,y\}$,  the dependent variable $h = h(x,y,
t)$ satisfies a fourth order partial derivative equation which is
of the form as follows:
\begin{eqnarray}
\frac{\partial h}{\partial t}=\nu \nabla^2 h-\kappa\nabla^4
h+\lambda {\mid \nabla h \mid} ^2,\;\;\;\;\;\;\;
\nabla=(\frac{\partial}{\partial x},\frac{\partial}{\partial y}).
\end{eqnarray}
where $h $  indicates the height of the interface.  The prefactor
$\nu$ is proportional to the surface tension coefficient and $\nu
\nabla^2 h$ is refered to as the surface tension term. The term
$\kappa\nabla^4 h$ is the result of surface diffusion which is
because of the curvature-induced chemical potential gradient. The
prefactor $\kappa$ is referred to as the surface diffusion term
and is proportional to the surface diffusion coefficient. The term
$\lambda {\mid \nabla h \mid} ^2$  represents the existence of
overhangs and vacancies during the deposition process (when
$\lambda
> 0$). The combination of the
$\nu \nabla^2 h$ and $\lambda {\mid \nabla h \mid} ^2$ terms
models the desorption effect of the deposited atoms \cite{[1]}.\\
Some important results concerning the K-S equation are already
known. For example, fractal tracer distributions of particles
advected in the velocity field of the K-S equation was found in
\cite{[4]}. Hong-Ji et al. in \cite{[1]} have studied the
evolution of (2+1)-dimensional surface morphology in the Kuramoto
- Sivashinsky K-S model by using the numerical simulation
approach.\\
The symmetry group method plays a fundamental role in the
analysis of differential equations.The theory of Lie symmetry
groups of differential equations was first developed by Sophus
Lie \cite{[7]} at the end of the nineteenth century, which was
called classical Lie method. Nowadays,  application of Lie
transformations group theory for constructing the solutions of
nonlinear partial differential equations (PDEs) can be regarded as
one of the most active fields of research in the theory of
nonlinear PDEs and applications.\\
The fact that symmetry reductions for many PDEs are unobtainable
by applying the classical symmetry method, motivated the creation
of several generalizations of the classical Lie group method for
symmetry reductions. The nonclassical symmetry method of reduction
was devised originally by Bluman and Cole in 1969 \cite{[8]}, to
find new exact solutions of the heat equation. The description of
the method is presented in \cite{[9],[10]}. Many authors have used
the nonclassical method to solve PDEs. In \cite{[11]} Clarkson and
Mansfield  have proposed an algorithm for calculating the
determining equations associated to the nonclassical method. A
new procedure for finding nonclassical symmetries has been
proposed by B$\hat{\mbox{i}}$l$\breve{\mbox {a}}$ and Niesen in \cite{[12]}.\\
$~~~$Classical and nonclassical symmetries of nonlinear PDEs may
be applied to reduce the number of independent variables of the
PDEs. Particularly, the PDEs can be reduced to ODES. The ODEs may
also have symmetries which enable  us to reduce the order of the
equation and we can integrate to obtain exact solutions. In
\cite{[5]}, R. Cimpoiasu et al. have discussed the classical and
nonclassical symmetries of the K-S model in the special case
$\nu=\kappa=\lambda=1$, but there are some problems because of
some computational mistakes. In this paper,  we will try  to
analyze the problem of the symmetries of the K-S equation in a
general case by considering all the three prefactors $\nu,
\kappa, \lambda$. Meanwhile, since this paper is the
generalization of \cite{[5]}, those computational mistakes exist
in \cite{[5]} will be corrected and some further facts are added.\\
$~~$ The structure of the present paper is as follows: In section
2, using the basic Lie symmetry method the most general Lie point
symmetry group of the K-S equation  is determined. In section 3,
some results yield from the structure of the Lie algebra of
symmetries are given.  Section 4 is devoted to obtaining the
one-parameter subgroups and the most general group-invariant
solutions of K-S equation. In section 5, we construct the optimal
system of one-dimensional subalgebras. Lie invariants and
similarity reduced equations corresponding to the infinitesimal
symmetries of equation (1.1) are obtained in section 6. In
section 7,  is devoted to the nonclassical symmetries of the  K-S
model, symmetries generated when a supplementary condition, the
invariance surface condition, is imposed. Some concluding remarks
are presented at the end of the paper.

%
\section{Classical symmetries of the K-S Equation}
In this section, the classical Lie symmetry method for the K-S
Equation has been performed. First, we recall the general
procedure for determining symmetries for an arbitrary system of
partial differential equations \cite{[1],[2]}. To begin, consider
a general system of partial differential equation containing $q$
dependent and $p$ independent variables as follows
\begin{eqnarray}
\Delta_\mu(x,u^{(n)})=0\ , \ \ \ \ \ \ \mu=1,...,r,
 \label{eq:7}
\end{eqnarray}
where $u^{(n)}$ represents all the derivatives of $u$ of all
orders from 0 to $n$. The one-parameter Lie group of
transformations
\begin{eqnarray}
\bar{x}^i=x^i+\varepsilon \xi^i (x,u)+O(\varepsilon^2),\quad
\bar{u}^\alpha=u^\alpha+\varepsilon \varphi^\alpha
(x,u)+O(\varepsilon^2),\qquad  i=1,...,p,\;\;\alpha=1,...,q.
\end{eqnarray}
where $\xi^i=\frac{\partial\bar{x}^i}{\partial
\varepsilon}|_{\varepsilon=0}$ and
$\varphi^\alpha=\frac{\partial\bar{u}^\alpha}{\partial
\varepsilon}|_{\varepsilon=0}$, are given. The action of the Lie
group can be recovered from that of its infinitesimal generators
acting on the space of independent and dependent variables.
Hence, we consider the following general vector field
\begin{eqnarray}
V=\sum_{i=1}^p \xi^i(x,u)\frac{\p}{\partial x^i}+\sum_{\alpha=1}^q
\varphi^{\alpha}(x,u) \frac{\partial}{\partial u^\alpha}
\end{eqnarray}
the characteristic of the vector field $V$ is given by the
function
\begin{eqnarray}
Q^\alpha(x,u^{(1)})=\varphi^\alpha(x,u)-\sum_{i=1}^p
\xi^i(x,u)\frac{\partial u^\alpha}{\partial x^i},\ \
\alpha=1,...,q.
\end{eqnarray}
Assume that the symmetry generator associated to (1.1) is given by
\begin{eqnarray}
V:=\xi^1(x,y,t,h)
\partial_x+\xi^2(x,y,t,h)\partial_y+\xi^3(x,y,t,h)\partial_t+\varphi(x,y,t,h)\partial_h\label{eq:3}
\end{eqnarray}
The fourth prolongation of $V$ is the vector field
\begin{eqnarray}
V^{(4)}=V+\varphi^x\partial_{h_x}+\varphi^t\partial_{h_y}+\varphi^t\partial_{h_t}+\varphi^{xx}\partial_{h_{xx}}+
\varphi^{xt}\partial_{h_{xt}}+ ...
+\varphi^{tttt}\partial_{h_{tttt}}\label{eq:5}
\end{eqnarray}
with coefficients given by
\begin{eqnarray}
\varphi^\iota=D_\iota Q+\xi^1 h_{x\iota}+\xi^2 h_{y\iota}+\xi^3
h_{t\iota},\quad\quad\quad
\varphi^{\iota\jmath}=D_{\imath}(D_\jmath Q)+\xi^1
h_{x\imath\jmath}+\xi^2 h_{y\imath\jmath}+\xi^3 h_{t\imath\jmath},
\end{eqnarray}
where $Q=\varphi-\xi^1 h_x-\xi^2 h_y-\xi^3 h_t$ is the
characteristic of the vector field $V$ given by (2.4) and $D_i$
represents total derivative and subscripts of $u$ are derivative
with respect to respective coordinates. $\imath$ and $\jmath$ in
the above could be $x$ and $t$ coordinates. By theorem (6.5) in
\cite{[]}, the invariance condition for the K-S is given by the
relation:
\begin{eqnarray}
\mathbf{V}^{(4)}[h_t+\nu
(h_{2x}+h_{2y})+\kappa(h_{4x}+2h_{(2x)(2y)}+h_{4y})-\lambda(h_x
^2+h_y ^2)]=0
\end{eqnarray}
 Hence, the invariance condition ()
is equivalent with the following equation:
\begin{eqnarray}\nonumber
&\varphi^t&+\nu(\varphi^{2x}+\varphi^{2y})+\kappa(\varphi^{4x}+2\varphi^{(2x)(2y)}+\varphi^{4y})-2\lambda(\varphi^x
h_x+\varphi^y u_y)=0
\end{eqnarray}
Substituting () into invariance condition (),
 we are left with a polynomial equation involving the
various derivatives of $h(x,y,t)$ whose coefficients are certain
derivatives of $\xi^1$, $\xi^2$, $\xi^3$ and $\varphi$. Since,
$\xi^1$, $\xi^2$, $\xi^3$ and $\varphi$ depend only on $x$, $y$,
$t$, $h$ we can equate the individual coefficients to zero,
leading to the complete set of determining equations:
\begin{eqnarray}\nonumber
&\xi^1 _h&=0,\quad\quad \xi^1 _x=0,\quad\quad \xi^1_y+\xi^2
_x=0,\quad\quad \xi^1 _{2t}=0\\\nonumber &\xi^2 _h&=0,\quad\quad
\xi^2 _y=0,\quad\quad \xi^2_{2t}=0,\quad\quad \xi^2
_{2x}=0\\\nonumber &\xi^3 _t&=0,\quad\quad \xi^3 _h=0,\quad\quad
\xi^3_{x}=0,\quad\quad \xi^3 _{y}=0,\quad\quad \xi^2 _{xt}=0
\\\nonumber&\varphi _t&=0,\quad\quad \varphi _h=0,\quad\quad
(2\lambda-2)\varphi_x-\xi^1
_t=0,\quad\quad(2\lambda-2)\varphi_y-\xi^2 _t=0
\end{eqnarray}
By solving this system of PDEs, we find that:
\begin{thm}
The Lie group of point symmetries of K-S (1.1) has a Lie algebra
generated by the vector fields
$\mathbf{V}=\xi\frac{\partial}{\partial
x}+\tau\frac{\partial}{\partial t}+\varphi
\frac{\partial}{\partial h}$, where
\begin{eqnarray}\nonumber
&\xi^1(x,y,t,h)&=c_4 y+c_2 t+c_3,\quad\quad \xi^2(x,y,t,h)=-c_4
x+c_5 t+c_6 ,\\\nonumber &\xi^3(x,y,t,h)&=c_1,\quad\quad
\varphi(x,y,t,h)=\big(\frac{x}{2\lambda}\big)c_2+\big(\frac{y}{2\lambda}\big)c_5+c_7.
\end{eqnarray}
and $c_i,\ i=1,...,7$ are arbitrary constants.
\end{thm}
\begin{cor}
Infinitesimal generators of every one parameter Lie group of point
symmetries of the K-S are:
\begin{eqnarray}\nonumber
&\mathbf{V}_1&=\frac{\partial}{\partial x},\;\;\;\;\;
\mathbf{V}_2=\frac{\partial}{\partial y},\;\;\;\;\;
\mathbf{V}_3=\frac{\partial}{\partial t},\;\;\;\;\;
\mathbf{V}_4=\frac{\partial}{\partial h},\;\;\;\;\;\\\nonumber
&\mathbf{V}_5&=y\frac{\partial}{\partial
x}-x\frac{\partial}{\partial y},\;\;\;\;\;
\mathbf{V}_6=t\frac{\partial}{\partial
x}+\frac{x}{2\lambda}\frac{\partial}{\partial h},\;\;\;\;\;
\mathbf{V}_7=t\frac{\partial}{\partial
y}+\frac{y}{2\lambda}\frac{\partial}{\partial h}
\end{eqnarray}
The commutator table of symmetry generators of the K-S  is given
in Table 1, where the entry in the $i^{\mbox{th}}$ row and
$j^{\mbox{th}}$ column is defined as $[V_i,V_j]=V_iV_j-V_j V_i,\
\ i,j=1,...,7.$ $\;\;\;\;\;\;\;\quad\quad\quad\Box$
\end{cor}
\begin{table}[h]
\centering{\caption{Commutation relations satisfied by
infinitesimal generators }}\label{table:1} \vspace{-0.35cm}
\begin{eqnarray*} \hspace{-0.75cm}\begin{array}{l |
l l l l l l l} \hline
                     & \hspace{0.6cm}V_1\hspace{1cm} & V_2\hspace{1cm}  & V_3\hspace{1cm} & V_4\hspace{1cm} & V_5\hspace{1cm} & V_6\hspace{1cm} & V_7\\ \hline
  V_1                  \hspace{0.7cm}  &\hspace{0.6cm} 0       & 0     &0         & 0 & {\bf V_2}& \hspace{-0.4cm}\big(\frac{1}{2\lambda}\big){\bf V_4} & 0\\
  V_2                       & \hspace{0.6cm} 0       & 0     &0 &0 & {\bf V_1}         &0 &\hspace{-0.4cm}\big(\frac{1}{2\lambda}\big){\bf V_4} \\
  V_3                        &\hspace{0.6cm} 0        & 0     & 0 &0  &0       &{\bf V_1}  &{\bf V_2}\\
  V_4                 & \hspace{0.6cm}0      & 0     &0         &0 &0 &0  &0\\
  V_5                 &\hspace{0.5cm} {\bf V_2}        &\hspace{-0.2cm} -{\bf V_1}    &0&0 &0  &\hspace{-0.19cm}{\bf V_7} & - {\bf V_6}\\
  V_6                 &\hspace{0.15cm}   \big(\frac{1}{2\lambda}\big){\bf V_4}       &0   &\hspace{-0.2cm}-{\bf V_1}       &0  &\hspace{-0.1cm}-{\bf V_7}&0 & 0 \\
  V_7                  & \hspace{0.6cm} 0       &\hspace{-0.7cm}  \big(\frac{-1}{2\lambda}\big){\bf V_4}     &\hspace{-0.2cm}-{\bf V_2} &0 &\hspace{-0.1cm} {\bf V_6} &0 &0\\
  \hline
 \end{array}\end{eqnarray*}
 \end{table}


\section{The Structure of the Lie algebra of Symmetries}
In this part, we determine the structure of symmetry Lie algebra
of the K-S equation.\\
$\goth{g}$ has no non-trivial {\it Levi decomposition} in the
form $\goth{g}=\goth{r}\ltimes \goth{g}_1$, because $\goth{g}$ has
not any non-trivial radical, i.e. if $\goth{r}$ be the radical of
$\goth{g}$, then $\goth{g}=\goth{r}$.





The Lie algebra $\goth{g}$ is solvable and non-semisimple. It is
solvable, because if
$\goth{g}^{(1)}=<V_i,[V_i,V_j]>=[\goth{g},\goth{g}]$, we have:
\begin{eqnarray*}
\goth{g}^{(1)}=[\goth{g},\goth{g}]=<-V_2, -
\frac{1}{2\lambda}V_4, V_1, V_7, -V_6>,\ \mbox{and}
\end{eqnarray*}
\begin{eqnarray*}
\goth{g}^{(2)}=[\goth{g}^{(1)},\goth{g}^{(1)}]=<\frac{1}{2\lambda}V_4>,
\end{eqnarray*}
so we have the following chain of ideals $\goth{g}^{(1)}\supset
\goth{g}^{(2)}\supset \{0\}$. Also, $\goth{g}$ is not semisimple,
because its killing form
\begin{eqnarray*}
\pmatrix{ 0&0&0&0&0&0&0\cr 0&0&0 &0&0&0&0\cr 0&0&0&0&0&0&0\cr
0&0&0&0&0 &0&0\cr 0&0&0&0&-4&0&0\cr 0&0&0&0&0&0&0 \cr
0&0&0&0&0&0&0}
\end{eqnarray*}
is degenerate.\\
Taking into account the table of commutators, $\goth{g}$ has two
abelian four and two dimensional subalgebras which are spanned by
$<V_1,V_2,V_3,V_4>$ and $<V_6,V_7>$, respectively, such that the
first one is an ideal in $\goth{g}$.

\section{Reduction of the Equation}
The equation (1.1) can be regarded as a submanifold of the jet
space $J^4(\Bbb{R}^3,\Bbb{R})$. So we can find the most general
group of invariant solutions of equation (1.1). To obtain the
group transformation which is generated by the infinitesimal
generators ${\bf V_i}=\xi^1 _i {\partial_x}+\xi^2 _i
{\partial_y}+\xi^3 _i
\partial _t+\varphi_i \partial _h$ for $i=1,...,7$, we need to solve the seven
systems of first order ordinary differential equations,
\begin{eqnarray}\nonumber
&\frac{d\bar{x}(s)}{ds}&=\xi^1_i(\bar{x}(s),\bar{y}(s),\bar{t}(s),\bar{h}(s)),\
\ \ \ \ \bar{x}(0)=x,\\
&\frac{d\bar{y}(s)}{ds}&=\xi^2_i(\bar{x}(s),\bar{y}(s),\bar{t}(s),\bar{h}(s)),\
\ \ \ \ \bar{y}(0)=y,\\
&\frac{d\bar{t}(s)}{ds}&=\xi^3_i(\bar{x}(s),\bar{y}(s),\bar{t}(s),\bar{h}(s)),\
\ \ \ \ \bar{t}(0)=t,\\\nonumber
&\frac{d\bar{u}(s)}{ds}&=\varphi_i(\bar{x}(s),\bar{y}(s),\bar{t}(s),\bar{u}(s)),\
\ \ \ \ \bar{u}(0)=u, \ \ \ \ \ \ \ \ i=1,...,7.
\end{eqnarray}
exponentiating the infinitesimal symmetries of equation (1.1), we
get the one parameter groups $G_k (x)$ generated by $v_k$ for
$k=1,...,7$.
\begin{thm}
The one-parameter groups $G_i(t):M\longrightarrow M$ generated by
the $V_i,\ i=1,...,7$, are given in the following table:
\begin{eqnarray*}
&G_1(s)&\ :\ (x,y,t,u)\longmapsto (x+s,y,t,u),\\ &G_2(s)&\ :\
(x,y,t,u)\longmapsto (x,y+s,t,u),\\\nonumber &G_3(s)&\ :\
(x,y,t,u)\longmapsto (x,y,t+s,u),\\
&G_4(s)&\ :\ (x,y,t,u)\longmapsto (x,y,t,u+s),\\\nonumber&G_5(s)&\
:\ (x,y,t,u)\longmapsto
 (x\cos(s)+y\sin(s),y\cos(s)-x\sin(s),t,u)\\\nonumber &G_6(s)&\ :\ (x,y,t,u)\longmapsto
 (x+st,y,t,\big(\frac{1}{4\lambda-4}\big)(ts^2+2xs+4u\lambda-4u)),\\\nonumber
 &G_7(s)&\ :\ (x,y,t,u)\longmapsto
(x,y+st,t,\big(\frac{1}{4\lambda-4}\big)(ts^2+2xs+4u\lambda-4u)).
\end{eqnarray*}
where entries give the transformed point
$\exp(tX_i)(x,y,t,u)=(\bar{x},\bar{y},\bar{t},\bar{u}).$
\end{thm}
These operators which generate invariance of the evolution
equation (2) corre- spond to the following transformations: U0
generate temporal translation, U1 and U4 Galilean boosts, U2
spatial dilatation and gauge transformation, U3 a gauge
transformation, U5 a spatial rotation, U6 and U7 spatial
translations.\\
 Recall that in general to each one parameter
subgroups of the full symmetry group of a system there will
correspond a family of solutions called invariant solutions.
Consequently, we can state the following theorem:
\begin{thm}
If $u=f(x,y,t)$ is a solution of equation (1.1), so are the
functions
\begin{eqnarray}
&G_1(s)&.f(x,y,t)=f(x+s,y,t),\\\nonumber
&G_2(s)&.f(x,y,t)=f(x,y+s,t),\\\nonumber
&G_3(s)&.f(x,y,t)=f(x,y,t+s),\\\nonumber
&G_4(s)&.f(x,y,t)=f(x,y,t)-s,\\\nonumber
&G_5(s)&.f(x,y,t)=f(x\cos(s)+y\sin(s),y\cos(s)-x\sin(s),t)
,\\\nonumber
&G_6(s)&.f(x,y,t)=f(x+st,y,t)-\big(\frac{2xs+ts^2}{4\lambda-4}\big),\\\nonumber
&G_7(s)&.f(x,y,t)=f(x,y+st,t)-\big(\frac{2ys+ts^2}{4\lambda-4}\big).\\\nonumber
\end{eqnarray}
\end{thm}
Thus, from above theorem we conclude that:
\begin{cor}
for the arbitrary combination $V =\sum_{i=1}^7 V_i\in\goth{g}$,
the K-S equation has the following solution
\begin{eqnarray}
u=e^{2\epsilon_7}f(xe^{\epsilon_7}+t\epsilon_5+\epsilon_1,ye^{\epsilon_7}\cos(\epsilon_6)+z\sin(\epsilon_6)+\epsilon_2,-yze^{\epsilon_7}\sin(\epsilon_6)
+z\cos(\epsilon_6)+\epsilon_3,te^{3\epsilon_7}+\epsilon_4)-\frac{\epsilon_5}{\alpha}
\end{eqnarray}
where $\epsilon_i$ are arbitrary real numbers.
\end{cor}
\section{Classification of Subalgebras for the K-S equation}
Let $\Delta$ be a system of differential equations with the
symmetry Lie group $G$.  Now,  $G$ operates on the set of the
solutions of $\Delta$ denoted by $S$. Let $s\cdot G$ be the orbit
of $s$, and $H$ be an $r-$dimensional subgroup of $G$. Hence, $H-$
invariant solutions $s\in S$ are characterized by equality
$s\cdot S=\{s\}$. If $h\in G$ is a transformation and $s\in S$,
then $h\cdot(s\cdot H)=(h\cdot s)\cdot (hHh^{-1})$. Consequently,
every $H-$invariant solution $s$ transforms into an $hHh^{-1}-$
invariant solution (Proposition 3.6 of \cite{[13]}).

Therefore, different invariant solutions are found from similar
subgroups of $G$. Thus, classification of $H-$invariant solutions
is reduced to the problem of classification of subgroups of $G$,
up to similarity. An optimal system of $r-$dimensional subgroups
of $G$ is a list of conjugacy inequivalent $r-$dimensional
subgroups of $G$ with the property that any other subgroup is
conjugate to precisely one subgroup in the list. Similarly, a
list of $r-$dimensional subalgebras forms an optimal system if
every $r-$dimensional subalgebra of $\goth g$ is equivalent to a
unique member of the list under some element of the adjoint
representation: $\tilde{\goth h}={\rm Ad}(g)\cdot{\goth h},\ g\in
G $.

Let $H$ and $\tilde{H}$ be connected, $r-$dimensional Lie
subgroups of the Lie group $G$ with corresponding Lie subalgebras
${\goth h}$ and $\tilde{\goth h}$ of the Lie algebra ${\goth g}$
of $G$. Then $\tilde{H}=gHg^{-1}$ are conjugate subgroups if and
only if $\tilde{\goth h}={\rm Ad}(g)\cdot{\goth h}$ are conjugate
subalgebras (Proposition 3.7 of \cite{[13]}). Thus, the problem of
finding an optimal system of subgroups is equivalent to that of
finding an optimal system of subalgebras, and so we concentrate on
it.
\subsection{Optimal system of one-dimensional subalgebras of the
K-S equation}
$~~$ There is clearly an infinite number of one-dimensional
subalgebras of the  K-S Lie algebra, $\goth{g}$, each of which
may be used to construct a special solutions or class of
solutions. So, it is impossible to use all the one dimensional
subalgebras of the  K-S to construct invariant solutions.
However, a well-known standard procedure \cite{[15]} allows us to
classify all the one-dimensional subalgebras into subsets of
conjugate subalgebras. This involves constructing the adjoint
representation group, which introduces a conjugate relation in
the set of all one-dimensional subalgebras.  In fact, for
one-dimensional subalgebras,the classification problem is
essentially the same as the problem of classifying the orbits of
the adjoint representation. If we take only one representative
from each family of equivalent subalgebras, an optimal set of
subalgebras is created. The corresponding set of invariant
solutions is then the minimal list from which we can get all
other invariant solutions of one-dimensional subalgebras simply
via transformations.\\
$~~$ Each $V_i,\ i=1,...,7$, of the basis symmetries generates an
adjoint representation (or interior automorphism)
$\mathrm{Ad}(\exp(\varepsilon V_i))$ defined by the Lie series
\begin{eqnarray}
\mathrm{Ad}(\exp(\varepsilon.V_i).V_j) =
V_j-s.[V_i,V_j]+\frac{\varepsilon^2}{2}.[V_i, [V_i,V_j]]-\cdots
\end{eqnarray}
where $[V_i,V_j]$ is the commutator for the Lie algebra, $s$ is a
parameter, and $i,j=1,\cdots,7$ (\cite{Olv1},page 199). In table 2
we give all the adjoint representations of the  K-S Lie group,
with the (i,j) the entry indicating $\mathrm{Ad}(\exp(\varepsilon
V_i))V_j$. Essentially, these adjoint representations simply
permute amongst "similar" one dimensional subalgebras. Hence,
they are used to identify similar one dimensional subalgebras.

We can expect to simplify a given arbitrary element,
\begin{eqnarray}
V=a_1 V_1+a_2V_2+a_3 V_3+a_4 V_4+a_5 V_5+a_6 V_6+a_7 V_7.
\end{eqnarray}
of the  K-S Lie algebra $\goth{g}$. Note that the elements of
$\goth{g}$ can be represented by vectors $a=(a_1,...,a_7)\in
{\Bbb R}^7$ since each of them can be written in the form (6.22)
for some constants $a_1,...,a_7$. Hence, the adjoint action can be
regarded as (in fact is) a group of linear transformations of the
vectors $(a_1,...,a_7)$.
\begin{table}[h]
\centering{\caption{Adjoint representation generated by the basis
symmetries of the  K-S Lie algebra }}\label{table:2}
\vspace{-0.35cm}
\begin{eqnarray*} \hspace{-0.75cm}\begin{array}{l |
l l l l l l l} \hline
  Ad                   & \hspace{1.5cm}V_1\hspace{3.7cm} & V_2\hspace{3.7cm}  & V_3\hspace{3.7cm} & V_4\\ \hline
  V_1                  \hspace{0.8cm}  & \hspace{1.4cm} V_1       & V_2    & V_3        &V_4 \\
  V_2                       & \hspace{0.7cm} \quad\quad V_1       & V_2    &V_3 &V_4  \\
  V_3                        &\hspace{0.7cm}\quad\quad V_1        & V_2    & V_3&V_4 \\
  V_4                 & \hspace{0.7cm}\quad \quad V_1     & V_2    &V_3         &V_4 \\
  V_5                 &\hspace{0.5cm} \cos(\varepsilon)V_1-\sin(\varepsilon)V_2        &\hspace{-0.8cm} \cos(\varepsilon)V_2+\sin(\varepsilon)V_1    &V_3&V_4  \\
  V_6                 &\hspace{0.6cm}   V_1+\big(\frac{\varepsilon}{2\lambda}\big){V_4}       &V_2   &\hspace{-0.4cm}V_3+\varepsilon V_1+\big(\frac{\varepsilon^2}{4\lambda-4}\big) V_4       &V_4 \\
  V_7                  & \hspace{1.2cm} V_1       &\hspace{-0.7cm} V_2+ \big(\frac{\varepsilon}{2\lambda}\big){ V_4}     &\hspace{-0.4cm}V_3+\varepsilon V_2+\big(\frac{\varepsilon^2}{4\lambda-4}\big)V_4  &V_4\\
  \hline
 \end{array}\end{eqnarray*}
 \end{table}
 \vspace{-0.7cm}
\begin{eqnarray*} \hspace{-0.75cm}\begin{array}{l |
l l l } \hline
  Ad                   & \hspace{1.5cm}V_5\hspace{3.7cm} & V_6\hspace{3.7cm}  & V_7\hspace{3.7cm} \\ \hline
  V_1                  \hspace{0.8cm}  &\hspace{1cm} V_5+\varepsilon V_2&\hspace{-0.4cm}  V_6-\big(\frac{\varepsilon}{2\lambda}\big){ V_4} & V_7   \\
  V_2                       & \hspace{1cm} V_5-\varepsilon V_1       &V_6 &\hspace{-0.4cm} V_7-\big(\frac{\varepsilon}{2\lambda}\big){ V_4}  \\
  V_3                        &\hspace{1.4cm}V_5       & \hspace{-0.3cm} V_6-\varepsilon V_1  &\hspace{-0.4cm}V_7-\varepsilon V_2 \\
  V_4                 & \hspace{1.4cm}  V_5 &V_6  &V_7 \\
  V_5                 & \hspace{1.4cm} V_5       &\hspace{-0.9cm} \cos(\varepsilon)V_6-\sin(\varepsilon)V_7    &\hspace{-0.5cm} \cos(\varepsilon)V_7+\sin(\varepsilon)V_6 \\
  V_6                 &\hspace{1cm}  V_5+\varepsilon V_7      &V_6   &V_7     \\
  V_7                  & \hspace{1cm}  V_5-\varepsilon V_6      &{ V_6}     &{ V_7}\\
  \hline
 \end{array}\end{eqnarray*}
Therefore, we can state the following theorem:
\begin{thm}
A one-dimensional optimal system of  the  K-S Lie algebra
$\goth{g}$ is given by
\begin{eqnarray}\nonumber
&(1)&:\ V_2+a V_6,\hspace{3.6cm}(4):\ aV_3+V_7,\\\nonumber
&(2)&:\ aV_3+bV_5, \hspace{3.4cm}(5):\ V_1+aV_3+bV_7,\\\nonumber
&(3)&:\ aV_3+V_6,\hspace{3.6cm} (6):\ aV_3+V_4+bV_5,\\\nonumber
\end{eqnarray}
where $a,b,c\in{\Bbb R}$ and $a\neq0$.
\end{thm}
{\it Proof:} $F^s_i:\goth{g}\to \goth{g}$ defined by
$V\mapsto\mathrm{Ad}(\exp(s_iV_i).V)$ is a linear map, for
$i=1,\cdots,7$. The matrix $M^s_i$ of $F^s_i$, $i=1,\cdots,7$,
with respect to basis $\{V_1,\cdots,V_7\}$ is
\begin{eqnarray}\nonumber
M_1 ^s=\pmatrix{1&0&0&0&0&0&0\cr 0&1&0 &0&0&0&0\cr
0&0&1&0&0&0&0\cr 0&0&0&1&0 &0&0\cr 0&s&0&0&1&0&0\cr
0&0&0&-s\zeta&0&1&0\cr 0&0&0&0&0&0&1} \ M_2^s=\pmatrix
 {1&0&0&0&0&0&0\cr 0&1&0
&0&0&0&0\cr 0&0&1&0&0&0&0\cr 0&0&0&1&0 &0&0\cr -s&0&0&0&1&0&0\cr
0&0&0&0&0&1&0 \cr 0&0&0&-s\zeta &0&0&1}\ M_3^s= \pmatrix{
1&0&0&0&0&0&0\cr 0&1&0 &0&0&0&0\cr 0&0&1&0&0&0&0\cr 0&0&0&1&0
&0&0\cr 0&0&0&0&1&0&0\cr -s&0&0&0&0&1&0\cr 0&-s&0&0&0&0&1}
\end{eqnarray}
\begin{eqnarray}\nonumber
\hspace{2cm} M_4^s=\pmatrix{ 1&0&0&0&0&0&0\cr 0&1&0 &0&0&0&0\cr
0&0&1&0&0&0&0\cr 0&0&0&1&0 &0&0\cr 0&0&0&0&1&0&0\cr 0&0&0&0&0&1&0
\cr 0&0&0&0&0&0&1}\quad\quad
 M_5 ^s=\pmatrix{C&-S&0&0&0&0&0\cr S&C&0&0&0&0&0\cr 0&0&1&0&0&0
&0\cr 0&0&0&1&0&0&0\cr 0&0&0&0&1&0&0 \cr 0&0&0&0&0&C&-S\cr
0&0&0&0&0&S&C}
\end{eqnarray}
\begin{eqnarray}\nonumber
\hspace{2cm} M_6 ^s=\pmatrix{1&0&0&s\zeta&0&0&0 \cr
0&1&0&0&0&0&0\cr s&0&1&s^2\frac{\zeta}{2}&0&0&0\cr 0&0&0&1&0&0&0
\cr 0&0&0&0&1&0&s\cr 0&0&0&0&0&1&0 \cr 0&0&0&0&0&0&1}\quad\quad
M_7 ^s=\pmatrix{ 1&0&0&0&0&0&0\cr 0&1&0 &s\zeta&0&0&0\cr
0&s&1&s^2\frac{\zeta}{2}&0&0&0\cr 0&0&0&1&0&0&0 \cr
0&0&0&0&1&-s&0\cr 0&0&0&0&0&1&0 \cr 0&0&0&0&0&0&1}
\end{eqnarray}
respectively, where $S=\sin s$, $C=\cos s$ and
$\zeta=\frac{1}{2\lambda}$. Let $V=\sum_{i=1}^7 a_i V_i$, then
\begin{eqnarray}
&& \hspace{-1mm} F^{s_7}\circ F^{s_6}_6\circ\cdots\circ
F^{s_1}_1\;:\;V\;\mapsto\; \\\nonumber
&&\hspace{4mm}\big(\cos(s_5)a_1-\sin(s_5)a_2+(\zeta
s_6\cos(s_5)-\zeta s_7\sin(s_5))a_4\big)V_1\\\nonumber
&&+\big(\sin(s_5)a_1+\cos(s_5)a_2+(\zeta s_7\cos(s_5)+\zeta s_6
\sin(s_5))a_4\big)V_2\\\nonumber && +\big(s_6 a_1+s_7
a_2+a_3+\frac{\zeta}{2}(s_7 ^2+s_6 ^2)a_4\big)V_3+a_4
V_4\\\nonumber
&&+\bigg(\big(-s_2\cos(s_5)+s_1\sin(s_5)\big)a_1+\big(s_2\sin(s_5)+s_1\cos(s_5)\big)a_2\\\nonumber
&&+\big((s_2\sin(s_5)+s_1\cos(s_5))\zeta
s_7+(-s_2\cos(s_5)+s_1\sin(s_5))\zeta s_6\big)a_4+a_5-s_7 a_6+s_6
a_7\bigg)V_5\\\nonumber
&&+\bigg(-s_3\cos(s_5)a_1+s_3\sin(s_5)a_2+\big(\zeta s_3 s_7
\sin(s_5)-\zeta s_3 s_6\cos(s_5)-\zeta
s_1\big)a_4+\cos(s_5)a_6-\sin(s_5)a_7\bigg)V_6\\\nonumber &&
+\bigg(-s_3\sin(s_5)a_1-s_3\cos(s_5)a_2+\big(-\zeta s_3
s_7\cos(s_5)-\zeta s_3 s_6\sin(s_5)-\zeta
s_2\big)a_4+\sin(s_5)a_6+\cos(s_5)a_7\bigg)V_7.
\end{eqnarray}
Now, we can simplify $V$ as follows:

If $a_4\neq 0$  we can make the coefficients of $V_1$, $V_2$ and
$V_6$ vanish by $F_1^{s_1}$, $F_2^{s_2}$, $F_6^{s_6}$ and
$F_7^{s_7}$. By setting $s_1=\frac{a_6}{\zeta a_4}$,
$s_2=\frac{a_7}{\zeta a_4}$, $s_6=-\frac{a_1}{\zeta a_4}$ and
$s_7=-\frac{a_2}{\zeta a_4}$, respectively. Scaling $V$ if
necessary, we can assume that $a_4=1$. So, $V$ is reduced to the
case (6).

If $a_4=0$ and $a_1\neq 0$ then we can make the coefficients of
$V_2$, $V_5$ and $V_6$ vanish by $F_5^{s_5}$, $F_3^{s_3}$ and
$F_2^{s_2}$. By setting $s_5=-\arctan(\frac{a_2}{ a_1})$,
$s_3=\frac{a_6}{ a_1}$ and $s_2=\frac{a_5}{a_1}$, respectively.
Scaling $V$ if necessary, we can assume that $a_1=1$. So, $V$ is
reduced to the case (5).

If $a_4=a_1=0$ and $a_2\neq 0$ then we can make the coefficients
of $V_3$, $V_7$ and $V_5$ vanish by $F_7^{s_7}$, $F_3^{s_3}$ and
$F_1^{s_1}$. By setting $s_7=-\frac{a_3}{ a_2}$, $s_3=\frac{a_7}{
a_2}$ and $s_1=-\frac{a_5}{a_2}$, respectively. Scaling $V$ if
necessary, we can assume that $a_2=1$. So, $V$ is reduced to the
case (1).

If $a_1=a_2=a_4=0$ and $a_6\neq 0$ then we can make the
coefficients of $V_7$ and $V_5$ vanish by $F_5^{s_5}$ and
$F_7^{s_7}$. By setting $s_5=-\arctan(\frac{a_7}{ a_6})$
 and $s_7=\frac{a_5}{a_6}$, respectively.
Scaling $V$ if necessary, we can assume that $a_6=1$. So, $V$ is
reduced to the case (3).

If $a_1=a_2=a_4=a_6=0$ and $a_7\neq 0$ then we can make the
coefficients of $V_5$  vanish by $F_6^{s_6}$. By setting
$s_6=-\frac{a_5}{a_7}$, respectively. Scaling $V$ if necessary, we
can assume that $a_7=1$. So, $V$ is reduced to the case (4).

If $a_1=a_2=a_4=a_6=a_7=0$ then  $V$ is reduced to the case (2).
$\;\;\;\;\;\;\;\Box$
\subsection{Two-dimensional optimal system}
The next step is  constructing the two-dimensional optimal system,
i.e., classification of two-dimensional subalgebras of $\goth{g}$.
This process is performed by selecting one of the vector fields
as stated in theorem (6). Let us consider $V_1$ (or
$V_i,i=2,3,4,5,6,7$). Corresponding to it, a vector field $V=a_1 V
_1+\cdots+a_7 V_7$, where $a_i$'s are smooth functions of
$(x,y,z,t)$ is chosen, so we must have
\begin{eqnarray}\label{eq:41}
[V_1,V]=\vartheta V_1+\varpi V.
\end{eqnarray}
Equation (\ref{eq:41}) leads us to the system
\begin{eqnarray}\label{eq:42}
C^i_{jk}\alpha_ja_k=\vartheta
a_i+\varpi\alpha_i\hspace{2cm}(i=1,\cdots,7).
\end{eqnarray}
The solutions of the system (\ref{eq:42}), give one of the
two-dimensional generator and the second generator is $V_1$ or,
$V_i,i=2,3,4,5,6,7$ if selected. After the construction of all
two-dimensional subalgebras, for every vector fields of theorem 6,
they need to be simplified by the action of adjoint matrices in
the manner analogous to the way of one-dimensional optimal system.
Thus the two-dimensional optimal system of $\goth{g}$ has three
classes of $\goth{g}$'s members combinations such as
\begin{eqnarray}\label{eq:43}
&<&\alpha_1 V _1+\alpha_2 V _2,\beta_1 V_3+\beta_2 V _4+\beta_3 V
_5>,\\\nonumber &<&\alpha_1 V _1+\alpha_2 V_3 ,\beta_1 V
_2+\beta_2 V _4+\beta_3 V_5>,\\\nonumber &<&V_1+\alpha V_3,
V_7>\;\;\;,\;\;\; <V_1+\alpha V_2,V_7>\\\nonumber &<&V_1+\alpha_1
V_6, \beta_1 V_4+\beta_2 V_5>\;\;\;\,;\;\;\; <V_4+\alpha_1
V_5+\alpha_2 V_6,V_1>.
\end{eqnarray}
where $\alpha_i,\  i=1,2$ and $\beta_j,\ j=1,..,3$ are real
numbers and $\alpha$ is a real nonzero constant. All of these
sub-algebras are abelian.
\subsection{Three-dimensional optimal system}
This system can be developed by the method of expansion of
two-dimensional optimal system. For this take any two-dimensional
subalgebras of (\ref{eq:43}), let us consider the first two vector
fields of (\ref{eq:43}), and call them $Y_1$ and $Y_2$, thus, we
have a subalgebra with basis $\{Y_1,Y_2\}$, find a vector field
$Y=a_1\textbf{V}_1+\cdots+a_7\textbf{V}_7$, where $a_i$'s are
smooth functions of $(x,y,z,t,u)$, such the triple $\{Y_1,Y_2,Y\}$
generates a basis of a three-dimensional algebra. For that it is
necessary an sufficient that the vector field $Y$ satisfies the
equations
\begin{eqnarray}\label{eq:44}
[Y_1,Y]=\vartheta_1Y+\varpi_1Y_1+\rho_1Y_2,\qquad[Y_2,Y]=\vartheta_2Y+\varpi_2Y_1+\rho_2Y_2,
\end{eqnarray}
and following from (\ref{eq:44}), we obtain the system
\begin{eqnarray}\label{eq:45}
C^i_{jk}\beta_r^ja_k=\vartheta_1a_i+\varpi_1\beta_r^i+\rho_1\beta_s^i,\qquad
C^i_{jk}\beta_s^ja_k=\vartheta_2a_i+\varpi_2\beta_r^i+\rho_2\beta_s^i.
\end{eqnarray}
The solutions of system (\ref{eq:45}) is linearly independent of
$\{Y_1,Y_2\}$ and give a three-dimensional subalgebra. This
process is used for the another two couple vector fields of
(\ref{eq:43}).

By performing the above procedure for all of the two couple vector
fields of (\ref{eq:43}), we conclude that $Y=\beta_1 Y_1+\beta_2
Y_2$.  By a suitable change of the base of $\goth g$, we can
assume that $Y=0$, so that $\goth g$ is not a $3-$dimensional
sub-algebra. Thus, we infer that:
\begin{cor}
The K-S equation Lie algebra $\goth{g}$ has no three-dimensional
Lie subalgebra.
\end{cor}
\section{Similarity Reduction of K-S equation} $~~~$ The
K-S equation (1.1) is expressed in the coordinates $(x,y,t,u)$, so
we ought to search for this equation's form in specific
coordinates in order to reduce it. Those coordinates will be
constructed by looking for independent invariants $(z,w,r)$
corresponding to the infinitesimal symmetry generator. Hence, by
applying the chain rule, the expression of the equation in the
new coordinate leads to the reduced equation.\\
$~~$We can now compute the invariants associated with the
symmetry operators. They can be obtained by integrating the
characteristic equations. For example for the  operator,
$H_5:=V_6= t\frac{\partial}{\partial
x}-\frac{1}{2\lambda}x\frac{\partial}{\partial h}$ this means:
\begin{eqnarray}
\frac{dx}{t}=\frac{dy}{0}=\frac{dt}{0}=\frac{-2\lambda dh}{x}
\end{eqnarray}
The corresponding invariants are as follows:
$z=y,\  w=t,\ r=h+\frac{x^2}{4\lambda t}$.\\
Taking into account the last invariant, we assume a similarity
solution of the form:
$h=f(z,w)-\frac{x^2}{4\lambda t}$
and we substitute it into (2) to determine the form of the
function $f(z,w)$: We obtain that $f(z,w)$ has to be a solution of
the following differential equation:
\begin{eqnarray}
f_w+\nu f_{zz}+kf_{4z}-\lambda f_z ^2-\frac{\nu}{2\lambda w}=0
\end{eqnarray}
Having determined the infinitesimals, the similarity solutions
$z_j$,  $w_j$ and $h_j$  are listed in Table I. In Table II we
list the reduced form of the K-S equation corresponding to
infinitesimal symmetries.
\begin{table}[h]
\centering{\caption{Lie Invariants and Similarity Solutions
}}\label{table:1} \vspace{-0.35cm}
\begin{eqnarray*} \hspace{-0.75cm}\begin{array}{l |
l l l l l l l} \hline
                 J    & \hspace{0.5cm}H_j\hspace{2.1cm} & z_j\hspace{2.1cm}  & w_j\hspace{2.1cm} & r_j\hspace{2.1cm} & h_j\\ \hline
 1                  \hspace{0.5cm}  &\hspace{0.5cm}  {\bf V_1}      & y    &t & h &f(z,w) \\
 2                     & \hspace{0.5cm} {\bf V_2}        & x    & t & h & f(z,w) \\
  3                        &\hspace{0.5cm} {\bf V_3}       & x   &y
  & h &
  f(z,w)\\
  4                 & \hspace{0.5cm}{\bf V_5}      &\hspace{-0.3cm} x^2+y^2    &t & h       &f(z,w) \\
  5                 &\hspace{0.5cm} {\bf V_6}        & y    &t& h+\frac{x^2}{4\lambda t} & f(z,w)-\frac{x^2}{4\lambda t} \\
   6                 &\hspace{0.5cm} {\bf V_7}        & x    &t& h+\frac{y^2}{4\lambda t}&f(z,w)-\frac{y^2}{4\lambda t} \\
  7                 &\hspace{0.2cm}   {\bf V_2+V_6}      &\hspace{-0.1cm}t   &y-\frac{x}{t}  &\hspace{-0.2cm}h+\frac{x^2}{4\lambda t}      &f(z,w)-\frac{x^2}{4\lambda t}   \\
  8                  & \hspace{0.2cm} {\bf V_3+V_5}        &\hspace{-0.3cm} x^2+y^2     &\hspace{-0.5cm} t-\arctan(\frac{x}{y})& h &f(z,w) \\
  9                  &\hspace{0.3cm} {\bf V_3+V_6}       & y  &\hspace{-0.3cm}-2x+t^2 & h+\frac{3tx-t^3}{6\lambda}       & f(z,w)-\frac{3tx-t^3}{6\lambda} \\
 10                   & \hspace{0.2cm} {\bf V_3+V_7}      & x   &\hspace{-0.1cm}-2y+t^2& h+\frac{3ty-t^3}{6\lambda} &f(z,w)-\frac{3ty-t^3}{6\lambda}  \\
  11                        &\hspace{0.2cm}  {\bf V_1+V_3+V_7}      & t-x     &\hspace{-0.2cm} y+\frac{1}{2}x^2-tx
  &h+\frac{x^3-3x^2 t+6xy}{12\lambda}
&f(z,w)-\frac{x^3-3x^2 t+6xy}{12\lambda}
    \\
  12                & \hspace{0.2cm}{\bf V_3+V_4+V_5}      & x^2+y^2    &\hspace{-0.8cm}t-\arctan(\frac{x}{y})  &\hspace{-0.3cm}h-\arctan(\frac{x}{y})    &f(z,w)+\arctan(\frac{x}{y}) \\
  \\
  \hline
 \end{array}\end{eqnarray*}
 \end{table}
\begin{table}[h]
\centering{\caption{Reduced equations corresponding to
infinitesimal symmetries }}\label{table:1} \vspace{-0.35cm}
\begin{eqnarray*} \hspace{-0.75cm}\begin{array}{l |
l l l l l l l} \hline
                 J    & \hspace{0.7cm}\mbox{Similarity Reduced Equations}\\ \hline
 1                  \hspace{0.7cm}  &\hspace{0.9cm}  f_{w}+\nu f_{zz}+k f_{4z}-\lambda f_z ^2=0\\

 2                     & \hspace{0.7cm} f_w+\nu f_{zz}+kf_{4z}-\lambda f_z^2=0 \\
   3                        &\hspace{0.7cm}\nu(f_{zz}+f_{ww})+k(f_{4z}+2f_{zzww}+f_{4w})-\lambda(f_{z}^2+f_w
   ^2)=0\\

    4                 & \hspace{0.7cm} f_w+4z\nu f_{zz}+4\nu f_z+16 k z^2 f_{4z}+64 z f_{zzz}+32 k f_{zz}-4\lambda zf_z^2=0 \\
  5                 &\hspace{0.7cm} f_w+\nu f_{zz}+kf_{4z}-\lambda f_z ^2-\frac{\nu}{2\lambda w}=0 \\
   6                 &\hspace{0.7cm} f_w+\nu f_{zz}+k f_{4z}-\lambda f_z^2-\frac{\nu}{2\lambda t}=0 \\
  7                 &\hspace{0.7cm}   2\lambda z^4 f_z+2\nu\lambda(z^2+z^4)f_{ww}+2k\lambda(z+1)^2f_{4w}-2\lambda^2(z^2+z^4)f_w ^2=0 \\
  8                  & \hspace{0.7cm}  (4k+\nu z)f_{ww}-\lambda z f_w ^2-4\lambda z^3 f_z ^2+64 k z^3 f_{3z}+16kz^4 f_{4z}+(4\nu z^3+32k z^2)f_{zz}
  +4\nu z^2 f_z+k f_{4w}+z^2 f_w+8kz^2 f_{zzww}=0\\
  9                 &\hspace{0.7cm} 16\nu\lambda f_{ww}-4\nu\lambda f_{zz}-64k\lambda f_{4w}-32k\lambda f_{zzww}-4k\lambda f_{4z}+16\lambda^2 f_{w}^2+4\lambda^2 f_z^2-w=0 \\
 10                   & \hspace{0.7cm} -4\nu\lambda f_{zz}-16\nu\lambda f_{ww}-4k\lambda f_{4z}-32k\lambda f_{zzww}-64k\lambda f_{4w}+4\lambda^2 f_z ^2+16\lambda^2f_w ^2-w=0 \\
  11                        &\hspace{0.7cm} 32k\lambda(z^2+z^4)f_{4w}+(96k\lambda z^2+32k\lambda)f_{3w}+(16z^2+48k\lambda+16\nu\lambda)f_{ww}+(96k\lambda z^2+32 k\lambda)
  f_{zzww}\\
 &\hspace{0.7cm}+64k\lambda(z+z^3)f_{zwww}+16\nu\lambda
 f_{zz}+(16\lambda-8\lambda w) f_z-16\lambda^2 (z^2+1)f_w
 ^2+(16\nu\lambda-16\lambda z w)f_w\\
&\hspace{0.7cm}-32 \lambda^2 z f_w f_z+192 k\lambda z f_{zww}+64
k\lambda z f_{zzzw}+32\nu\lambda z f_{zw}+8\nu z+4 w^2=0
    \\
  12                & \hspace{0.7cm} (4\nu z^3+32k z^2)f_{zz}+(z^2+2\lambda z)f_w+16 k z^4 f_{4z}+(\nu z +4k)f_{ww}+kf_{4w}-4\lambda z^3 f_z ^2+64k z^3 f_{3z}\\
  &\hspace{0.7cm}     +4\nu z^2 f_z+8kz^2 f_{zzww}-\lambda z f_w ^2 -\lambda z=0   \\
  \\
  \hline
 \end{array}\end{eqnarray*}
 \end{table}

\section{Nonclassical Symmetries of  K-S Equation}
 In this section,  we will apply the so called nonclassical symmetry
method  \cite{[8]}. Beside the classical symmetries, the
nonclassical symmetry method can be used to find some other
solutions for a system of PDEs and ODEs. The nonclassical
symmetry method has become the focus of a lot of research and
many applications to physically important partial differential
equations as in \cite{[9],[10],[11],[12]}. Here, we follow the
method used by  {\sc cai guoliang} et al, for obtaining the
non-classical symmetries of the Burgers-Fisher equation based on
compatibility of evolution equations \cite{[17]}. For the
non-classical method, we must add the invariance surface
condition to the given equation, and then apply the classical
symmetry method. This can also be conveniently written as:
\begin{eqnarray}
V^{(4)}\Delta_1|_{\Delta_1=0,\Delta_2=0}=0,\label{eq:15}
\end{eqnarray}
where $V$ is defined in (2.5) and $\Delta_1$ and $\Delta_2$ are
given as:
\begin{eqnarray*}
\Delta_1:=h_t+h_{2x}+h_{2y}+h_{4x}+2h_{(2x)(2y)}+h_{4y}-(h_x
^2+h_y ^2) ,\hspace{1cm} \Delta_2:=\varphi-\xi^1 h_x-\xi^2
h_y-\xi^3h_t
\end{eqnarray*}
Without loss of generality we choose $\xi^3=1$. In this case using
$\Delta_2$ we have:
\begin{eqnarray}
 h_t=\varphi-\xi^1 h_x-\xi^2 h_y.
\end{eqnarray}
First, total differentiation $D_t$ of the equation gives
\begin{eqnarray*}
 D_t(
h_t)=D_t(-\nu h _{2x}-\nu h_{2y}-k h_{4x} -2k h_{(2x)(2y)}-k
h_{4y}+\lambda h_x ^2+\lambda h_y ^2)
\end{eqnarray*}
\vspace{-0.9cm}
\begin{eqnarray*}
 D_t(\varphi-\xi^1 u_x-\xi^2
u_y)&=&-\nu h _{2xt}-\nu h_{2yt}-k h_{4xt} -2k h_{(2x)(2yt)}-k
h_{4yt}+2\lambda h_x  h_{xt}+2\lambda
h_y h_{yt}\\
&=&-\nu D_{xx}(\varphi-\xi^1 u_x-\xi^2 u_y)-\nu
D_{yy}(\varphi-\xi^1
u_x-\xi^2 u_y)-k D_{xxxx}(\varphi-\xi^1 u_x-\xi^2 u_y)\\
 &-&k
D_{yyyy}(\varphi-\xi^1 u_x-\xi^2 u_y)-2k D_{xxyy}(\varphi-\xi^1
u_x-\xi^2 u_y)+2\lambda h_x D_{x}(\varphi-\xi^1 u_x-\xi^2 u_y)\\
&+&2\lambda h_y D_{y}(\varphi-\xi^1 u_x-\xi^2 u_y).
\end{eqnarray*}
Substituting $\xi^1 u_{xt}$  and $\xi^2 u_{yt}$ to both sides, we
can get
\begin{eqnarray*}
\varphi^t&=&-\nu\varphi^{xx}-\nu\varphi^{yy}-k\varphi^{xxxx}-2k\varphi^{xxyy}-k\varphi^{yyyy}+2\lambda
h_x \varphi^x+2\lambda h_y \varphi^y+\xi^1 u_{xt}+\xi^2 u_{yt}\\
&+&\xi^1 u_{xxx}+\xi^2 u_{xxy}+\xi^1 u_{yyx}+\xi^2 u_{yyy}+k\xi^1
u_{xxxxx}+k\xi^2 u_{xxxxy}+2k\xi^1
u_{xxxyy}\\&+&\xi^2u_{xxyyy}-2\lambda \xi^1 u_x u_{xx}-2\lambda
\xi^2 u_x u_{xy}-2\lambda u_y \xi^1 u_{xy}-2\lambda u_y \xi^2
u_{yy}
\end{eqnarray*}
By virtue of
\begin{eqnarray*}
 D_x(
h_t)=D_x(-\nu h _{2x}-\nu h_{2y}-k h_{4x} -2k h_{(2x)(2y)}-k
h_{4y}+\lambda h_x ^2+\lambda h_y ^2)\\
 D_y(
h_t)=D_y(-\nu h _{2x}-\nu h_{2y}-k h_{4x} -2k h_{(2x)(2y)}-k
h_{4y}+\lambda h_x ^2+\lambda h_y ^2)
\end{eqnarray*}
gives
\begin{eqnarray*}
&h_{xt}&=-\nu h _{xxx}-\nu h_{2yx}-k h_{xxxxx} -2k h_{(2x)(2yx)}-k
h_{4yx}+2\lambda h_x  h_{xx}+2\lambda h_y h_{yx}\\
&h_{yt}&=-\nu h _{2xy}-\nu h_{yyy}-k h_{4xy} -2k h_{(2x)(yyy)}-k
h_{yyyyy}+2\lambda h_x  h_{xy}+2\lambda h_y h_{yy}
\end{eqnarray*}
so it gives the governing equation
\begin{eqnarray*}
\varphi^t+\varphi^{2x}+\varphi^{2y}+\varphi^{4x}+2\varphi^{(2x)(2y)}+\varphi^{4y}-2\lambda(\varphi^x
u_x+\varphi^y u_y)=0
\end{eqnarray*}
where $\varphi^t$, $\varphi^x$ are given by
\begin{eqnarray*}
\varphi^t&=&D_t(\varphi-\xi^1 u_x-\xi^2 u_y-\xi^3 u_t)+\xi^1 u_{xt}+\xi^2 u_{yt}+\xi^3 u_{tt}=D_t(\varphi-\xi^1 u_x-\xi^2 u_y)+\xi^1 u_{xt}+\xi^2 u_{yt},\\
\varphi^x&=&D_x(\varphi-\xi^1 u_x-\xi^2 u_y-\xi^3 u_t)+\xi^1 u_{xx}+\xi^2 u_{xy}+\xi^3 u_{xt}=D_x(\varphi-\xi^1 u_x-\xi^2 u_y)+\xi^1 u_{xx}+\xi^2 u_{xy},\\
\varphi^y&=&D_y(\varphi-\xi^1 u_x-\xi^2 u_y-\xi^3 u_t)+\xi^1
u_{xy}+\xi^2 u_{yy}+\xi^3 u_{ty}=D_t(\varphi-\xi^1 u_x-\xi^2
u_y)+\xi^1 u_{xy}+\xi^2 u_{yy},\\
\varphi^{xx}&=&D_{xx}(\varphi-\xi^1 u_x-\xi^2 u_y-\xi^3 u_t)+\xi^1
u_{xxx}+\xi^2 u_{xxy}+\xi^3 u_{xxt}=D_{xx}(\varphi-\xi^1 u_x-\xi^2
u_y)+\xi^1 u_{xxx}+\xi^2 u_{xxy},\\
\varphi^{yy}&=&D_{yy}(\varphi-\xi^1 u_x-\xi^2 u_y-\xi^3 u_t)+\xi^1
u_{xyy}+\xi^2 u_{yyy}+\xi^3 u_{tyy}=D_{yy}(\varphi-\xi^1 u_x-\xi^2
u_y)+\xi^1 u_{xyy}+\xi^2 u_{yyy},\\
\varphi^{xxxx}&=&D_{xxxx}(\varphi-\xi^1 u_x-\xi^2 u_y-\xi^3
u_t)+\xi^1
u_{xxxxx}+\xi^2 u_{yxxxx}+\xi^3 u_{txxxx}\\
&=&D_{xxxx}(\varphi-\xi^1
u_x-\xi^2 u_y)+\xi^1 u_{xxxxx}+\xi^2 u_{yxxxx},\\
\varphi^{xxyy}&=&D_{xxyy}(\varphi-\xi^1 u_x-\xi^2 u_y-\xi^3
u_t)+\xi^1
u_{xxxyy}+\xi^2 u_{xxyyy}+\xi^3 u_{xxyyt}\\
&=&D_{xxxx}(\varphi-\xi^1
u_x-\xi^2 u_y)+\xi^1 u_{xxxxx}+\xi^2 u_{yxxxx},\\
\varphi^{yyyy}&=&D_{yyyy}(\varphi-\xi^1 u_x-\xi^2 u_y-\xi^3
u_t)+\xi^1
u_{xyyyy}+\xi^2 u_{yyyyy}+\xi^3 u_{tyyyy}\\
&=&D_{yyyy}(\varphi-\xi^1 u_x-\xi^2 u_y)+\xi^1 u_{xyyyy}+\xi^2
u_{yyyyy}.
\end{eqnarray*}
Substituting them into the governing equation, we can get the
 determining equations for the symmetries of the K-S equation.
  Substituting $\xi^3=1$ into the
determining equations we obtain the determining equations of the
nonclassical symmetries of the original equation (1.1). Solving
the system obtained by this procedure, the only solutions we found
were exactly the solution  obtained through the classical
symmetry approach (theorem 1). This means that no supplementary
symmetries, of non-classical type, are specific for our model.


\section*{Conclusion}
$~~$In this paper by applying the criterion of invariance of the
equation under the infinitesimal prolonged infinitesimal
generators, we find the most general Lie point symmetries group of
the  K-S equation. Also, we have constructed the optimal system
of one-dimensional subalgebras of  K-S equation. The latter,
creats the preliminary classification of group invariant
solutions. The Lie invariants and similarity reduced equations
corresponding to infinitesimal symmetries are obtained. By
applying the nonclassical symmetry method  for the  K-S model we
concluded that the analyzed model do not admit supplementary,
nonclassical type symmetries. Using this procedure, the classical
Lie operators only were generated.

\end{document}